\documentclass[12pt]{article}

\usepackage{amsfonts}
\usepackage{amsmath}
\usepackage{amssymb}
\usepackage{amsthm}
\usepackage[line,matrix,arrow]{xy}
\usepackage{graphics}

\newtheorem{thm}{Theorem}[section]

\theoremstyle{plain}

\newtheorem{lemma}[thm]{Lemma}

\newtheorem{gevolg}[thm]{Corollary}
  \newtheorem{opques}[thm]{Open Problem}

\numberwithin{equation}{section}

\newcommand{\re}{\stackrel{\rightarrow}{e}}
\newcommand{\el}{\stackrel{\leftarrow}{e}}

\theoremstyle{remark}

\newtheorem{rem}[thm]{Remark}

\begin{document}

\title{A modified version of frozen percolation on the binary tree}

\author{R. Brouwer \\
{\small CWI, Amsterdam } \\
{\footnotesize email: Rachel.Brouwer@cwi.nl}
}
\date{}

\maketitle

\begin{abstract}
We consider the following, intuitively described process: at time zero, all sites of a binary tree are at rest. Each site becomes activated at a random uniform $[0,1]$ time, independent of the other sites. As soon as a site is in an infinite cluster of activated sites, this cluster of activated sites freezes. The main question is whether a process like this exists. Aldous~\cite{Aldous} proved that this is the case for a slightly different version of frozen percolation. In this paper we construct a process that fits the intuitive description and discuss some properties.  
\end{abstract}

\begin{section}{Introduction}
Let $G$ be an arbitrary graph and let $N \geq 2$ be an integer. We attach a random uniform $[0,1]$ number $U_i$ to each site $i$ of $G$, independently. 
Informally, the frozen percolation model can be described as follows. 
At time $0$ all sites are at rest, or coloured white. At time $U_i$, the site $i$ becomes active, or coloured green. If a site is in a cluster of green sites that has size at least $N$, this cluster {\em freezes}, or becomes red, instantaneously. The term frozen percolation was firstly introduced by Aldous~\cite{Aldous} in a slightly different setting\footnote{Aldous freezes infinite clusters. Readers familiar with Aldous' work should note that in his description, the boundary of a frozen cluster is by definition at rest at the time a cluster freezes, and stays at rest forever. In our description however, at the time when a frozen cluster arises, its boundary can be either at rest (and allowed to get activated and even frozen later on), or already frozen. }. We discuss his results briefly in the next section \ref{froztree}.

We use both the terminology of freezing and of the colours in what follows, whichever seems more appropriate. So, how does our model evolve in time? Some sites get activated and never freeze, others will freeze at the moment that they are activated and some after they were activated. It follows from the above description that at time $1$ all sites are green or red. 

It should be clear that the frozen percolation model exists for finite graphs $G$ and any $N$. In particular, when $N$ is larger than the number of sites in $G$, each site turns green at some point, and stays green forever. As soon as we take the underlying graph infinite, existence of the model is not immediately clear. 

The existence of the frozen percolation model is relatively easy on $\mathbb{Z}$. For any time $t < 1$, there will be infinitely many sites $j$ on the left and right half lines that have $U_j \geq t$. This `breaks up' the line in finite pieces; the state of any site at time $t$ then depends on finitely many variables only and this is sufficient to show existence.  More or less the same reasoning can be used to show that the process can be constructed on $\mathbb{Z}^d$, for fixed $N$.

Now we consider infinite graphs, it is tempting to take $N$ equal to infinity, i.e.~only infinite clusters are allowed to freeze. On $\mathbb{Z}$, this does not lead to interesting behaviour: all sites in $\mathbb{Z}$ stay green up to time $1$ and at time $1$, every site becomes red. On a (regular) binary tree, the process does exist in the slightly different form presented by Aldous~\cite{Aldous}. Benjamini and Schramm have shown (but not published) that on $\mathbb{Z}^2$, the frozen percolation process where infinite clusters freeze (in Aldous' setting), does not exist. See also~\cite{BergToth}.  
In this paper we construct a 
frozen percolation model where infinite clusters freeze on the tree, using the dynamics described above.  
It is a priori not clear that the frozen percolation process where infinite clusters freeze exists on the regular binary tree. In a slightly different setting Aldous~\cite{Aldous} proved that a form of this process indeed exists and we will discuss his result briefly in the next subsection. Later we discuss the existence of the process on the binary tree for our version of frozen percolation.

\end{section} 

\begin{section}{Description of Aldous' result}
Informally, the model considered by Aldous is as follows. Assign a uniform, independent $0$-$1$ variable $U_e$ to each edge $e$ of the regular binary tree. Let $\mathcal{A}_0 = \emptyset$. For each edge $e$, at time $t = U_e$ set $\mathcal{A}_{t} = \mathcal{A}_{t^-} \cup \{e\}$ if each end-vertex of $e$ is in a finite cluster of $A_{t^-}$; otherwise set $\mathcal{A}_{t} = \mathcal{A}_{t^-}$. That is, the boundaries of infinite clusters will never join the process. In the final configuration, there are infinite clusters, finite clusters and boundary edges in between. Aldous shows that this process exists by first `guessing' what the distribution function of the time that an edge joins the process should be on the directed tree. Heuristic arguments suggest that this distribution function should be 
\begin{equation}
G(t) = \left\{
\begin{array}{ll}
0 & t \in [0,1/2) \\
1 - \frac{1}{2t} & t \in [1/2,1]\\
\frac{1}{2} & 1 < t < \infty \\
1 & t = \infty.
\end{array} \right.
\label{AldFt}
\end{equation} 
Then \eqref{AldFt} is used to construct a process on the undirected tree and finally Aldous shows that the model indeed meets with its intuitive description. We follow this line of reasoning below.   

The model proposed by Aldous is meant as a model for polymerisation. In this model a polymer is made up of molecular units, where each unit is capable of forming three bonds. Before the critical time, there are only finite polymers (the sol) and later on infinite polymers (the gel) are observed. Aldous' model displays the behaviour observed by chemical physicists: beyond the gel point, the number of small polymers decreases but their average size retains a constant value. Formally, this means that {\em at any time $t> t_c$} (where $t_c$ corresponds to the critical time of ordinary percolation on the tree), finite clusters of edges that have joined the process have the same distribution as critical percolation clusters. See Proposition 11 of~\cite{Aldous}. It is a perfect example of self-organised criticality (although Aldous does not mention SOC).  

A priori, we see no reason why the boundary of infinite clusters should not join the process. Computations turn out to be harder if they do and we do not get such nice results as Aldous did. Nevertheless, it seems worth to try to construct the frozen percolation process where boundaries of infinite clusters are allowed to join the process. 
From now on, we refer to frozen percolation using our dynamics as modified frozen percolation to distinguish between the two versions.
It is interesting to see whether the self-organised critical behaviour observed in Aldous' model is `a coincidence', in the sense that a small perturbation of the dynamics causes the critical behaviour to disappear. In the same spirit (see also Remark \ref{edge}), we chose to study the site version of the modified frozen percolation model. We will see in Section \ref{crittree} that, although the perfect correspondence to ordinary critical percolation is lost, the modified frozen percolation process still behaves critically.  
\end{section}

\begin{section}{Heuristics and computational arguments on the directed tree}

Consider the directed rooted binary tree $\mathcal{T}$, whose root $O$ has degree $2$ and all other vertices have degree $3$. Let each edge be directed away from the root. To each site, we assign a random uniform $[0,1]$ number $U_i$ independently. The dynamics is as before, but we now freeze infinite directed rays: at time $0$ all sites are at rest (coloured white). Each site $i$ becomes activated (green) at time $t=U_i$, but as soon as $i$ is in an infinite (directed) active cluster, this cluster freezes (becomes red) instantaneously. 
In this section, we argue what the behaviour of the directed modified frozen percolation process should be like. Later on, we use these heuristic arguments to formally construct the modified frozen percolation process.
 
Let $1$ and $2$ be the children of $O$. We use $\mathcal{T}[i]$, $i = 1,2$ to denote the subtree that has site $i$ as its root, see Figure \ref{bintree2}.
\begin{figure}[h]
\centerline{\scalebox{0.5}{\includegraphics{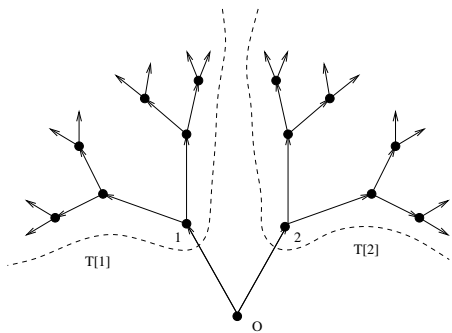}}}
\caption{The rooted binary tree $\mathcal{T}$ and its subtrees $\mathcal{T}[1]$ and $\mathcal{T}[2]$.}
\label{bintree2}
\end{figure}
If the root $O$ is activated, it freezes as soon as one of its children freezes. Suppose $Y$ is the time that site $O$ freezes. If $O$ never freezes, we take $Y$ equal to infinity. Since $\mathcal{T}[1]$ and $\mathcal{T}[2]$ are isomorphic to $\mathcal{T}$, we require that the time $Y_{1}$ ($Y_{2}$) that site $1$ ($2$) freezes in $\mathcal{T}[1]$ ($\mathcal{T}[2]$, respectively) is distributed like $Y$. 
Furthermore, a site cannot freeze before it is in an infinite cluster of activated sites so $Y$ should be at least $1/2$, by ordinary percolation results. 
To formalise these statements, define $I$ to be $[1/2,1] \cup \{\infty\}$. 
Let $\Phi(x,y,z)$ be the following function on $I \times I \times [0,1]$.
\begin{eqnarray}
\Phi(x,y,z) = \begin{cases}
x \text{        if } x \geq z \\
y \text{        if } x < z \leq y \\
\infty \text{     otherwise.}
\end{cases}
\label{defPhi}
\end{eqnarray}     
The heuristic arguments above require that  
\begin{equation}
Y \stackrel{d}{=} \Phi(\min\{Y_{1},Y_{ 2}\}, \max\{Y_{1},Y_{2}\},U_O). 
\label{recursive}
\end{equation}
\noindent 
The following lemma shows which distribution functions satisfy equality \eqref{recursive}.  
\begin{lemma} Let $F$ be a (possibly defective) probability distribution function, of an $I$-valued r.v. with the following additional properties: $F$ is continuous and differentiable with strictly positive derivative on $[1/2,1]$ and has $F(1/2) = 0$. 
Let $(Y_1,Y_2,U)$ be independent random variables, $Y_1,Y_2$ each having probability distribution function $F$ and $U$ having the uniform distribution on $[0,1]$. Then, 
\begin{equation}
\Phi(\min(Y_1,Y_2), \max(Y_1,Y_2),U) \text{ again has probability distribution function } F, 
\nonumber
\end{equation}
if and only if 
\begin{equation}
F(t) = \left\{ \begin{array}{ll} 
0  & t \leq 1/2 \\
\ln(2t) & 1/2 < t \leq 1 \\
\ln(2) &  1 < t < \infty \\
1 & t = \infty.
\end{array}
\right.
\label{defF}
\end{equation}
\label{consistency}
\end{lemma}
\begin{proof}
The equation $F(t) = \mathcal{P}(\Phi(\min(Y_1,Y_2),\max(Y_1,Y_2),U) \leq t)$ is equivalent to 
\begin{equation}
F(t) = \int^t_{0} \mathcal{P}(t \geq \min(Y_1,Y_2) > s)ds\ + \int^t_{0} \mathcal{P}(\min(Y_1,Y_2) < s < \max(Y_1,Y_2) \leq t) ds.
\nonumber 
\end{equation}
Since both $Y_1$ and $Y_2$ have probability distribution function $F$ this is equivalent to 
\begin{eqnarray}
F(t)&=& \int^t_{0} \left[2F(t) -F(t)^2 - 2F(s) + F(s)^2\right] ds + 2\int^t_{0} (F(t) - F(s))F(s) ds \nonumber \\
&=& t[2F(t) - F(t)^2] + 2F(t) \int^t_{0} F(s)ds - \int^t_{0} 2F(s) + F(s)^2 ds.
\label{check}
\end{eqnarray}
Differentiating with respect to $t$ we obtain 
\begin{equation}
\frac{d F(t)}{d t} = \frac{d F(t)}{d t} [2t(1-F(t))] + 2\frac{d F(t)}{d t} \int^t_{0} F(s) ds.
\label{2sol}
\end{equation}

%

Since $F$ has positive derivative on $[1/2,1]$, we may divide \eqref{2sol} by $\frac{d F(t)}{d t}$. The remaining integral equation has $F(t) = \ln(t) + C$, for some constant $C$, as its solution. 
Finally, we use $F(1/2) = 0$ to obtain $C = \ln(2)$. 

The reversed implication follows from straightforward calculation. 
\end{proof}
Note that the solution above is defective, i.e.~has mass in infinity. We can interpret this distributions as the probability that a vertex freezes before time $t$, as we have seen above. 
 
\medskip

When we release the restriction that $F$ must be differentiable with positive derivative on $[1/2,1]$ in the lemma above, we obtain more solutions. 
One of these is 
\begin{equation}
F(t) = \left\{ \begin{array}{ll}
0 & t < \infty \\
1 & t = \infty, 
\end{array}
\right.
\label{zero-sol}
\end{equation}
We have the following heuristic argument to rule out the zero solution \eqref{zero-sol} as a candidate. We approximate the directed process by considering processes on finite directed trees. As before, an independent uniform $[0,1]$ variable $U_i$ is attached to each site $i \in \mathcal{T}$. Let $\mathcal{P}$ denote the measure governing the $U$-variables in $\mathcal{T}$. The root of the tree $O$ is said to be at level $0$ and a site is at level (depth) $n$ if its distance to the root is $n$. For each $n$, the dynamics is as follows. Let $\mathcal{T}(n)$ denote the first $n+1$ levels of the tree. At time $0$ all sites are at rest. Each site becomes activated at time $U_i$ and as soon as there is a directed path of activated sites to the leaves of the tree (i.e. level $n$) this path freezes. For fixed $n$ this process exists and as we did on $\mathbb{Z}$, we can define sets evolving in time to formalize the dynamics. 
Let $\mathcal{W}_n(0) := V_{\mathcal{T}(n)}$, all vertices of $\mathcal{T}(n)$, $\mathcal{G}(0) := \emptyset$ and $\mathcal{R}(0) := \emptyset$. For $t > 0$, 
\begin{eqnarray}
\mathcal{W}_n(t) &:=& \{ i \in \mathcal{T}(n) : U_i > t \}, \nonumber \\
\mathcal{G}_n(t) &:=& \{ i \in \mathcal{T}(n) : U_i \leq t \text{ and for all paths } \pi = i,i_2,\ldots,i_k \text{ with $i_k$ on level $n$},\nonumber \\ && \exists j \in \{2,\ldots,k\} \text{ with } U_j > t \text{ or } U_{i_k} < \max \{U_i,U_{i_2},\ldots,U_{i_{k-1}}\})\nonumber \\
\mathcal{R}_n(t) &:=& \{ i \in \mathcal{T}(n) : U_i \leq t \text{ and } i \notin \mathcal{G}(t))
\nonumber \\
\end{eqnarray}
Analogously to the previous section we define 
\begin{equation}
F_n(t) := \mathcal{P}(O \in \mathcal{R}_n(t)). 
\nonumber
\end{equation}
\begin{lemma}
For all $t \in [1/2,1]$, we have 
\begin{equation}
\limsup_{n\rightarrow \infty} F_n(t) \geq 1 - \frac{1}{2t}.
\label{>theta/2}
\end{equation}
\label{FN}
\end{lemma}
Note that the lower bound in \eqref{>theta/2} is exactly half the ordinary percolation function on the binary tree. Further, it is equal to the distribution $G$ \eqref{AldFt} on the directed tree that drops out in Aldous´ model, corresponding to the distribution $F$ \eqref{defF} in our case.  The proof of this lemma uses the same ideas as the proof of Lemma 4.5 of~\cite{forestfires} and is in fact simpler so we do not present it here.  

At this point, it is not clear whether $F_n(t)$ converges as $n$ tends to infinity. Even if it does, it is not clear that its distribution converges to the infinite process we described above. 
Nevertheless, Lemma \ref{FN} suggests that we can rule out the solution \eqref{zero-sol}. In view of the heuristic arguments, we propose the distribution $F$ \eqref{defF}.  
as a candidate for the directed modified frozen percolation process. 
From now on, let $F$ denote the distribution function in \eqref{defF}.

So far, our arguments (apart from Lemmas \ref{consistency} and \ref{FN}) are non-rigorous. We have inferred what the distribution of the time that a site freezes should be, if the modified frozen percolation process exists on the directed binary tree. However, there are still many questions. Does the modified frozen percolation process exist on the directed binary tree? If so, can we use the directed process to construct the modified frozen percolation process on the Bethe lattice $\mathcal{B}$? On the Bethe lattice, where each site has degree three, we can consider each site as being attached to the roots of three directed subtrees. It is not clear however, whether the directed processes exist simultaneously: does there exists a unique translation invariant law of the states of the sites at time $t, 0 \leq t \leq 1$, for {\em all} directions simultaneously? Finally, if the modified frozen percolation process can be constructed on $\mathcal{B}$, does it display critical behaviour? 

In the following section we construct the modified frozen percolation model on $\mathcal{B}$. We follow Aldous~\cite{Aldous} but we need to make several non-trivial adjustments.  
\label{heuristics}
\end{section}




\begin{section}{The frozen percolation process on $\mathcal{B}$} 
%
From now on we consider the Bethe lattice $\mathcal{B}$ where each site has degree $3$. There is no natural sense of direction as there was on the rooted directed tree, so that we cannot apply the results from the previous section directly. We solve this by directing the adjoined edges of a site $i$ outward. The children of $i$, called $i_1$,$i_2$ and $i_3$ here, form the roots of directed subtrees $\mathcal{T}[i_1]$, $\mathcal{T}[i_2]$ and $\mathcal{T}[i_3]$, where the direction of the edges is inherited from the edges leaving from $i$, see Figure \ref{dir3}.
\begin{figure}[h]
\centerline{\scalebox{0.5}{\includegraphics{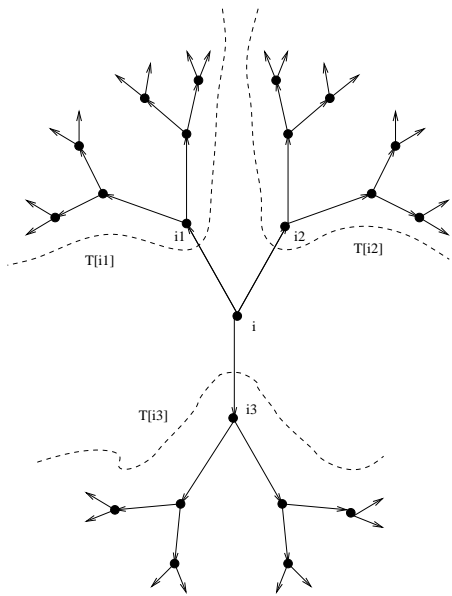}}}
\caption{The site $i$, and the three subtrees $\mathcal{T}[i_1]$, $\mathcal{T}[i_2]$ and $\mathcal{T}[i_3]$ leaving from its children.}
\label{dir3}
\end{figure}

We define $Y_{i \rightarrow i_1}$ to be the time that $i_1$ freezes in the directed subtree $\mathcal{T}[i_1]$. We use the subscript $i \rightarrow i_1$ to indicate that 
we consider ${i_1}$ to be the root of the directed subtree leaving from $i_1$, where the direction of the edges is inherited from the edge $i \rightarrow i_1$. The variables $Y_{i \rightarrow i_2}$ and $Y_{i \rightarrow i_3}$ should be likewise. In what follows, if we write $Y_{i \rightarrow j}$ it is implicit that $i$ is a neighbour of $j$. 

The following lemma shows that we can define the variables $Y_{i \rightarrow j}$ for all $i \in \mathcal{B}$ and all $j \in \partial \{i\}$, in a consistent manner. 
Recall that $\partial\{i\}$ denotes the set of children of $i$. 
\begin{lemma} There exists a joint law for $((U_{i},Y_{i\rightarrow j}): i \in \mathcal{B}, j \in \partial \{i\})$ which is invariant under automorphisms of the tree and such that for each $i \in {\mathcal{B}}$ and each $j \in \partial \{i\}$ we have 
\renewcommand{\theenumi}{(\roman{enumi})}
\begin{enumerate}
\item{\begin{equation} Y_{j \rightarrow i} \text{ has distribution function } F. \nonumber \end{equation}}
\item{
\begin{eqnarray} 
Y_{j \rightarrow i} = \Phi(\min(Y_{i \rightarrow k},Y_{i \rightarrow l}), \max(Y_{i \rightarrow k},Y_{i \rightarrow l}),U_i) \text{  a.s., }
\label{defy}
\end{eqnarray}
where $\{k,l\} = \partial \{i\} \backslash \{j\}$. 
}
\item{For each finite connected set $S \subset \mathcal{B}$, the variables $(Y_{i \rightarrow j} : i \in S, j \notin S)$ are independent of each other and independent of the collection $(U_i : i \in S)$. 
\label{ind}}
\end{enumerate}

\label{existencelemma}
\end{lemma}
\begin{proof}
Let $O$ be the root and let $n \geq 1$. Define $V_{\leq n}$ to be the set of sites at distance at most $n$ from $O$. Let $V_{n}$ be the set of sites at distance $n$. Take $(U_i : i \in V_{\leq n})$ independent of each other. Then take $(Y_{j \rightarrow i} :  j \in V_{n}, i \in V_{n+1})$ independent of the above mentioned $U_i$'s and independent of each other, with each $Y_{j \rightarrow i}$ having probability distribution $F$. We apply equation \eqref{defy} successively to define $Y_{i \rightarrow j}$ for all sites $i,j \in V_{\leq n}$ and by Lemma \ref{consistency} all these have law $F$.
As $n$ increases, the joint laws are consistent and hence we can apply the Kolmogorov consistency theorem to obtain a joint law for  $((U_{i},Y_{i\rightarrow j}): i \in \mathcal{B}, j \in \partial \{i\})$. Automorphism invariance of this law is straightforward. Part \ref{ind} of the lemma easily follows from the construction. 
%
\end{proof}
\begin{gevolg}
\begin{equation}
\mathcal{P}(\exists i,k \in \mathcal{B}, l \in \partial\{k\} \text{ s.t. } U_i = Y_{k \rightarrow l}) = 0.
\nonumber 
\end{equation}
\label{notequal}
\end{gevolg}
\begin{proof}
Since the number of sites in $\mathcal{B}$ is countable it is sufficient to show that for fixed $i,k$ and $l$ with $l \in \partial \{k\}$, $U_i \neq Y_{k \rightarrow l}$ almost surely. By automorphism invariance, we may assume that $i = O$. We take $n$ so large that the site $k$ is in $V_{\leq n}$. By the construction in the proof of Lemma \ref{existencelemma}, $Y_{k \rightarrow l}$ must be equal to one of the $Y$-values in $(Y_{v \rightarrow w}:  v \in V_{n}, w \in V_{n+1})$ or infinity. Since $n$ is fixed, we can modify on a null set so that $U_i$ is different from the values $(Y_{v \rightarrow w} : v \in V_{n}, w \in V_{n+1})$ almost surely. It follows that $U_i \neq Y_{k \rightarrow l}$, almost surely.  
\end{proof}



The process on $\mathcal{B}$ can now be described by the following: 
We define for each $i \in \mathcal{B}$,  
\begin{eqnarray}
Z_i &:=& \min\{ Y_{j \rightarrow i} : j \in \partial \{i\} \} \nonumber \\  
    & =& \min\{ Y_{i \rightarrow j} : j \in \partial\{i\}, Y_{i \rightarrow j} \geq U_i \},
\label{defZ}
\end{eqnarray}
where we take $Z_i = \infty$ if the minimum does not exist. The equality follows from \eqref{defy}. The last expression seems more complicated but we use it later on, since we want to make use of the fact that $U_i$ is independent of $Y_{i \rightarrow j}$. 

We define 
$\mathcal{W}(0) := V_{\mathcal{B}}, \mathcal{G}(0) := \emptyset$ and $\mathcal{R}(0) := \emptyset$. 
These evolve in time as follows: for $t \in [0,1]$, 
\begin{eqnarray}
\mathcal{W}(t) &:=&  \{ i \in \mathcal{B} : U_i > t \}, \nonumber \\
\mathcal{G}(t) &:=&  \{ i \in \mathcal{B} : U_i \leq t, Z_i > t \},  \nonumber \\
\mathcal{R}(t) &:=&  \{ i \in \mathcal{B} : Z_i \leq t \}. 
\label{defg(t)}
\end{eqnarray}

One should think of the set $\mathcal{W}(t)$ as the set of white (at rest) sites;  of $\mathcal{G}(t)$ as the set of green (activated but in a finite cluster) sites and $\mathcal{R}(t)$ as the set of sites that are red (frozen infinite clusters) at time $t$. We use this terminology from now on. 

From the definition immediately follows that for all sites $i$,   with $Z_i < \infty$, 
\begin{eqnarray}
&& i \in \mathcal{W}(t) \text{ for } t \in [0,U_i), \nonumber \\
&& i \in \mathcal{G}(t) \text{ for } t \in [U_i,Z_i), \nonumber \\
&& i \in \mathcal{R}(t) \text{ for } t \in [Z_i,1]. \nonumber
\end{eqnarray}
If $Z_i = \infty$, clearly 
\begin{eqnarray}
&& i \in \mathcal{W}(t) \text{ for } t \in [0,U_i), \nonumber \\
&& i \in \mathcal{G}(t) \text{ for } t \in [U_i,1]. \nonumber
\end{eqnarray}
 
This description gives us the following: 
\begin{lemma} For every vertex $v$ and for every pair of neighbouring vertices $v \sim w$, 
\renewcommand{\theenumi}{(\roman{enumi})}
\begin{enumerate}
\item{The probability that a vertex is eventually not frozen equals $$\mathcal{P}(v \in \mathcal{G}(1)) = \frac{3}{2}\ln(2)^2 - \frac{1}{2} \approx 0.22, $$
so that the probability that a vertex is eventually frozen equals $$\mathcal{P}(v \in \mathcal{R}(1)) = 1 - \mathcal{P}(v \in \mathcal{G}(1)) \approx 0.77. $$
\label{vfreexe}
}
\item{The probability that both $v$ and $w$ are eventually frozen but in different frozen infinite clusters equals $$\mathcal{P}(v,w \in \mathcal{R}(1), Z_v \neq Z_w) = 3\ln(2) - 2 \approx 0.079.$$}
\end{enumerate}
\label{explcom}
\end{lemma}
\begin{proof} 
The lemma follows easily from the definition of $Z_i$'s \eqref{defZ} and $Y_{\cdot \rightarrow \cdot}$'s \eqref{defy} and the distribution $F$ \eqref{defF}. For example, part \ref{vfreexe}: 
\begin{eqnarray}
\mathcal{P}(v \in \mathcal{G}(1)) &=& \mathcal{P}(Z_v = \infty) \nonumber \\
&=& \mathcal{P}(\forall j \in \partial\{v\} : Y_{v \rightarrow j} < U_i \text{ {\em or} } Y_{v \rightarrow j} = \infty) \nonumber \\
&=& \int_{1/2}^1[\ln(2s) + 1 - \ln(2)]^3 ds + \frac{1}{2}(1-\ln(2))^3 = \frac{3}{2}\ln(2)^2 - \frac{1}{2}. \nonumber 
\end{eqnarray}
The other statement is proved similarly.   
\end{proof}
  
In contrast to the finite tree where it was clear that the sets $\mathcal{W}(t), \mathcal{G}(t)$ and $\mathcal{R}(t)$ fit their intuitive description, we need to prove this in some detail here. So far, we have constructed {\em a} process, but it is not clear that it acts like the informally described modified frozen percolation process. 
The following lemmas show that only infinite clusters join $\mathcal{R}(\cdot)$ and that there are no infinite green clusters in $\mathcal{G}(\cdot)$, almost surely.  

\begin{lemma}
Let $S_s(i)$ denote the cluster of $i$ in $\mathcal{G}(s)$, considered as a set of sites. Almost surely, $\forall i \in \mathcal{B}$ with $Z_i < \infty$,  
\renewcommand{\theenumi}{(\roman{enumi})}
\begin{enumerate}
\item{$S_{Z_i^-}(i) \subseteq \mathcal{R}(Z_i)$ \label{totaljoin}}
\item{$|S_{Z_i^-}(i)| = \infty$\label{infjoin}}
\end{enumerate}
\label{redok}
\end{lemma}
\begin{proof} 
Since the number of sites in $\mathcal{B}$ is countable, it is sufficient to show that a.s. \ref{totaljoin} and \ref{infjoin} hold for some fixed $i$. Choose an arbitrary site $i$, with $Z_i < \infty$. 
By definition of $Z_i$ \eqref{defZ}, there exists a site $j \in \partial\{i\}$ such that $Y_{i \rightarrow j} \geq U_i$ and $Z_i = Y_{i \rightarrow j}$. 
Using \eqref{defy} repeatedly there exists an infinite self-avoiding path $i =: j_0,j =: j_1, j_2, \ldots $ such that 
\begin{eqnarray}
&& Z_i = Y_{i \rightarrow j} = Y_{j \rightarrow j_2} = Y_{j_2 \rightarrow j_3} = \cdots , 
\label{Yarray}
\end{eqnarray}
and 
\begin{eqnarray}
U_{j_k} < Z_i \text{ for } k \geq 0, \text{ almost surely.} \label{Uj<t} 
\end{eqnarray}
Note that we have used Corollary \ref{notequal} to obtain a strict inequality in \eqref{Uj<t}. 
In the directed subtree where the direction is given by $i \rightarrow j$, each site $j_{k}$ has only two children. One of them is $j_{k+1}$; let $v_{k+1}$ be the other child, for $k > 0$. 
By \eqref{defy} we have for all $k > 0$, 
\begin{equation}
Y_{j_{k} \rightarrow v_{k+1}} > Z_i \text{ or } Y_{j_{k} \rightarrow v_{k+1}} < U_{j_k}. 
\label{fex}
\end{equation} 
We need to show that the path $j_1,j_2,\ldots$ constructed above has $Z_{j_k} = Y_{j_{k} \rightarrow j_{k+1}} (= Z_i)$ for all $k>0$. 
By definition 
\begin{equation}
Z_{j_k} = \min \{ Y_{j_k \rightarrow w} : w \in \{j_{k+1}, j_{k-1}, v_{k+1}\}, Y_{j_k \rightarrow w} \geq U_{j_k} \}.  \nonumber  
\end{equation} 
By \eqref{Yarray} and \eqref{fex}, the minimum cannot be achieved taking $w$ equal to $v_{k+1}$. So suppose that $U_{j_k} \leq Y_{j_k \rightarrow j_{k-1}} < Z_i$. Applying \eqref{defy} and \eqref{fex} repeatedly we obtain 
\begin{equation}
Y_{j_k \rightarrow j_{k-1}} = Y_{j_{k-1} \rightarrow j_{k-2}} = \cdots = Y_{j \rightarrow i} = Y_{i \rightarrow w} < Z_i, 
\label{===<t}
\end{equation}
for some $w \in \partial\{i\}, w \neq j$. From \eqref{defy} also follows that $Y_{i \rightarrow w} \geq U_i$, which implies $Z_i \leq Y_{i \rightarrow w} < Z_i$, a contradiction. 
This shows that $Z_i = Z_{j} = Z_{j_2} = \cdots$ which together with \eqref{Uj<t} shows that an infinite path contained in $S_{Z_i^-}(i)$ joins $\mathcal{R}(Z_i)$ at time $Z_i$. This proves part \ref{infjoin} of the lemma.

Now suppose that we have a site $v_k$ in $S_{Z_i^-}(i)$ that is adjacent to the infinite path $j_0,j_1,\ldots$ at the site $j_{k-1}$ for some $k>0$. Because $v_k \in S_{Z_i^-}(i)$, we have that $U_{v_k} < Z_i$ and further that $Z_{v_k} \geq Z_i$.  
We need to show that equality holds here. But this follows easily from the above, since $Y_{v_k \rightarrow j_{k-1}} = Y_{j_{k-1}\rightarrow {j_k}} =  Z_{j_{k-1}} = Z_i$ by the construction of the infinite path above, and hence $Y_{v_k \rightarrow j_{k-1}} \geq U_{v_k}$ and 
\begin{equation}
Z_i \leq Z_{v_{k}} \leq Y_{v_k \rightarrow j_{k-1}} = Z_{j_{k-1}} = Z_i. 
\nonumber 
\end{equation}
We repeat this argument `working to the outside of $S_{Z_i^-}(i)$' so that we obtain that $S_{Z_i^-}(i) \subseteq \mathcal{R}(Z_i)$, which proves part \ref{totaljoin} of the lemma.  
\end{proof} 

\begin{lemma}
Almost surely, there is no infinite component in $\mathcal{G}(t)$ for $t \in [0,1)$.  
\label{upto1}
\end{lemma}
\begin{proof}
Suppose there exist $t \in [0,1)$ and $i \in \mathcal{B}$ such that 
$$S_{t}(i) \subseteq \mathcal{G}(t) \text{ and } |S_{t}(i)| = \infty.$$ 
Then $Z_i > t$, by definition. We will show that almost surely, all $j \in S_{t}(i)$ have $Z_j = Z_i$, so that if an infinite green cluster exists at some time-point, it will exist in a time-interval of positive length. 
Suppose that there is a $k \in S_{t}(i)$ with $Z_k \in (t,Z_i)$. This means that $S_t(i) = S_t(k) \subseteq S_{Z_k^-}(k) \subseteq \mathcal{R}(Z_k)$ almost surely, by Lemma \ref{redok}. So $i \in \mathcal{R}(Z_k)$ and this contradicts the fact that $Z_i > Z_k$.  
We conclude that almost surely, all $j \in S_t(i)$, $Z_j \geq Z_i$ so that $S_{t}(i) \subseteq \mathcal{G}(Z_i^-)$.  

From the above argument follows that if $\mathcal{G}(t)$ contains an infinite component, then there exist two rationals $t_1 < t_2$ such that $\mathcal{G}(s)$ contains an infinite component for all $s \in [t_1,t_2]$. Since the number of such pairs $t_1,t_2$ is countable, it is sufficient to show that for fixed $t_1<  t_2$,  
\begin{equation}
\mathcal{P}( \exists \text{ an infinite component which, } \forall s \in [t_1,t_2], \text{ is in } \mathcal{G}(s)) = 0. 
\label{toQ} 
\end{equation}
For $t_1 \in [0,1/2]$ the above easily follows from ordinary percolation results. 

Now suppose $t_1 \in [1/2,1)$ and let $t_2 \in (t_1,1]$. Then fix a site $i_0$ and a path $\pi = i_0, i_1, \ldots, i_{n-1}$ of length $n$. Add direction to the subtree leaving from $i_1$ consistent with $i_0 \rightarrow i_1$.    
Let for $j \in [2,n-1]$, $v_{j}$ denote the child of $i_{j-1}$ that is not on $\pi$. 
Consider   
\begin{eqnarray} 
&&\hspace{-30pt} \mathcal{P}(\forall s \in [t_1,t_2]\ \forall i \in \pi : i \in \mathcal{G}(s)) \nonumber \\
&=&\mathcal{P}(\forall j \in [0,n-1] : U_{i_j} \leq t_1 \text{ and } Z_{i_j} \geq t_2)  \nonumber \\
&\leq& \mathcal{P}(\forall j \in [1,n-2]: U_{i_j} \leq t_1 \text{ and either } Y_{i_j \rightarrow v_{j+1}} \geq t_2 \text{ or } Y_{i_j \rightarrow v_{j+1}} < U_{i_j} )  
\nonumber \\
&=& \bigl[\mathcal{P}(U_{i_1} \leq t_1 \text{ and either } Y_{i_1 \rightarrow v_2} \geq t_2 \text { or } Y_{i_1 \rightarrow v_2} < U_{i_1})\bigr]^{n-2}. \label{inbetween}
\end{eqnarray}
The inequality follows from \eqref{defZ} and the equality from automorphism invariance and the independence property, see Lemma \ref{existencelemma}, part \ref{ind}. Now using the knowledge about the distribution of $Y_{i_1 \rightarrow v_2}$ (recall that it is equal to $F$, \eqref{defF}) the last term in \eqref{inbetween} equals 
\begin{equation}
\bigl[ \int_{1/2}^{t_1}\ln(2s)ds + t_1[1-\ln(2t_2)] \bigr]^{n-2} 
= \bigl[ 1/2 + t_1 \ln (t_1/t_2)\bigr]^{n-2}
\nonumber
\end{equation}
Summing over all $3\cdot 2^n$  possible paths leaving from $i_0$ we obtain 
\begin{eqnarray} 
&&\hspace{-20pt} \mathcal{P}(\exists \text{ path $\pi$ of length $n$ leaving from } i_0 \text{ such that } \forall j \text{ on } \pi:  U_{j} \leq t_1 \text{ and } Z_{j} \geq t_2) \nonumber \\
 &\leq& 12 \bigl( 1 + 2t_1 \ln(t_1/t_2) \bigr)^{n-2}. 
\label{n-2}
\end{eqnarray}
Note that $t_2 > t_1$ so that $\ln(t_1/t_2) < 0$. Applying \eqref{n-2}, we get 
\begin{eqnarray}
&& \hspace{-30pt}\mathcal{P}( \exists \text{ infinite path $\pi$ from $i_0$ which, }\forall s \in [t_1,t_2], \text{ is contained in }\mathcal{G}(s) ) 
\nonumber \\
 &\leq& \lim_{n \rightarrow \infty} 12 \bigl( 1 + 2t_1 \ln(t_1/t_2) \bigr)^{n-2} = 0. 
\nonumber
\end{eqnarray}
We sum over all $i_0$ to show \eqref{toQ} for $t_1 \in (1/2,1)$. This finishes the proof of the lemma.   
\end{proof}



Summarising, we have constructed a process on the binary tree with the following properties. At time $0$, all sites are at rest. Each site $i$ becomes activated at time $U_i$, where $U_i$ is a uniform $U[0,1]$ random variable (chosen independently for each $i$). As soon as a site is in an infinite activated cluster at some time $t < 1$, it freezes. 
In the following section we will see that the modified frozen percolation process displays some features of critical behaviour. 
\end{section}
 
\begin{section}{Critical behaviour}
Aldous 
shows 
for his dynamics 
that finite non-empty clusters (which are the analogue of green clusters in our terminology) are distributed as ordinary critical percolation clusters conditional on being non-empty. Remarkably this holds for {\em all} $t \geq 1/2$. His frozen percolation process behaves like a critical system at all times (after the critical time), so that it is a nice example of self-organised criticality. 

In our case, the distribution of green clusters is different: suppose that $S$ is an arbitrary finite connected set of sites such that $|S| > 1$; let $\partial S$ denote its boundary. If $v$ is a neighbour of $w$ we write $v \sim w$. 
\begin{lemma}
$S \subseteq \mathcal{G}(t)$  if and only if
\renewcommand{\theenumi}{(\alph{enumi})} 
\begin{enumerate}
\item{$\forall w \in S : U_w \leq t$ and \label{iffset1}}
\item{$\forall v \notin S \text{ such that } v \sim s \text{ for some }s \in S, \text{ we have either } Y_{s \rightarrow v} < U_s \\
\text{ or } Y_{s \rightarrow v} > t.$ 
\label{iffset2}}
\end{enumerate}
\label{g(t)}
\end{lemma}
\begin{proof}
One of the implications is trivial. Suppose $S \subseteq \mathcal{G}(t)$. Then \ref{iffset1} and \ref{iffset2} hold by definition of $\mathcal{G}(t)$ \eqref{defg(t)}. The reversed implication is only slightly more difficult. Suppose that \ref{iffset1} and \ref{iffset2} hold. Then, `working from the outside in', we can recursively find all $Y_{w_1 \rightarrow w_2}$ for all $w_1 \in S$, $w_2 \in \partial\{w_1\}$ using \eqref{defy}. It is easy to see that 
all $Y$-values thus obtained are in $(t,1]$ or equal to $\infty$, so that $S \subseteq \mathcal{G}(t)$. 
\end{proof}
   
If we want to compute the probability that $S \subseteq \mathcal{G}(t)$, i.e.~by Lemma \ref{g(t)} the joint probability of \ref{iffset1} and \ref{iffset2} above, the actual geometry of the set $S$ may be of importance. Not surprisingly, for $t \leq 1/2$ we find 
$$ \mathcal{P}(S \subseteq \mathcal{G}(t)) = t^{|S|}.$$ 
Now suppose $t > 1/2$ and consider \ref{iffset1} and \ref{iffset2}. 
Every $w \in S$ that has two neighbours $v_1, v_2$ in $\partial S$, contributes a factor 
\begin{equation}
\int_{1/2}^t \ln(2s)^2 ds + 2(1-\ln(2t))\int_{1/2}^t \ln(2s) ds + t(1-\ln(2t))^2 = t - \ln(2t).  
\nonumber 
\end{equation} 
Note that the contribution of such $w$ to the joint probability of \ref{iffset1} and \ref{iffset2} only depends on $U_w, Y_{w \rightarrow v_1}$ and $Y_{w \rightarrow v_2}$, by part \ref{ind} of Lemma \ref{existencelemma}. Intuitively, this factor can be explained by considering $w$ to be activated, but not frozen in the directed subtree where $w$ is the root, and $v_1$ and $v_2$ are its children. 

Further, every $w$ that has only one neighbour $v$ in $\partial S$ contributes a factor 
$$\int_{1/2}^t\ln(2s) ds + t(1-\ln(2t)) = 1/2.$$ 
The contribution of such $w$ to the joint probability of \ref{iffset1} and \ref{iffset2} only depends on $U_w$ and $Y_{w \rightarrow v}$. The occurrence of the factor `$1/2$' is quite surprising. We have no intuitive explanation, but it is the reason for the critical behaviour observed (see for example equation \eqref{1/2} and Lemma \ref{expcrit}) later on. 
Finally, every $w \in S$ that has no neighbours in $\partial S$ contributes a factor $t$. 
Define
\begin{eqnarray}
n(0) := \# \{w \in S: \# \{v \in \partial S, v \sim w\} = 0 \}, \nonumber \\
n(1) := \# \{w \in S: \# \{v \in \partial S, v \sim w\} = 1 \}, \nonumber \\
n(2) := \# \{w \in S: \# \{v \in \partial S, v \sim w\} = 2 \}. \nonumber 
\end{eqnarray}
Note that for any connected set $S$ we have $n(0) + 2 = n(2)$.  
For any connected set $S$ such that $|S| > 1$, 
\begin{equation}
\mathcal{P}(S \subseteq  \mathcal{G}(t)) = (t-\ln(2t))^{n(2)}\bigl(\frac{1}{2}\bigr)^{n(1)}t^{n(0)}. 
\label{compclusterdistr}
\end{equation}


From this we can see that the distribution of green clusters is not the same as distribution of ordinary percolation clusters (for any parameter value), because not only the size but also the geometry of a cluster plays a role. Nevertheless, by \eqref{compclusterdistr} we have for all $t \geq 1/2$, and all sites $v$ and $w$ 
\begin{equation}
\mathcal{P}(\exists \text{ green path from } v \text{ to } w \text{ at time }t) = \left(\frac{1 - \ln(2)}{2}\right)^2 \left(\frac{1}{2}\right)^{\pi(v,w)}, 
\label{1/2}
\end{equation} 
where $\pi(v,w)$ is the number of sites in the path from $v$ to $w$. This is the most remarkable similarity to ordinary critical percolation; there the probability that the path from $v$ to $w$ is contained in an open cluster equals $(1/2)^{\pi(v,w)}$.  

Arbitrary green clusters are also critical in some sense: 
\begin{lemma}
Let $E$ denote the expectation with respect to the law of the frozen percolation process. Fix a site $O$ and let $\mathcal{G}_O(t)$ denote the green cluster of $O$ at time $t$. We consider $V_{n}$, the set of sites at distance $n$ from $O$. Then there exist positive constants $A_1$ and $A_2$ such that for all $n$, 
\begin{equation}
A_1 \leq E| \mathcal{G}_O(t) \cap V_n | \leq A_2. 
\end{equation}
\label{expcrit}
\end{lemma} 
\begin{proof}
\begin{eqnarray}
E|\mathcal{G}_O(t) \cap V_{n} | &=& \sum_{v \in V_{n}} \mathcal{P}(v \in \mathcal{G}_O(t)) \nonumber \\ 
&=& \sum_{v \in V_n} \mathcal{P}(\text{the unique path from $O$ to $v$ is in $\mathcal{G}(t)$})  \nonumber \\ 
&=&  3 \cdot 2^{n-1} \bigl(\frac{1}{2}\bigr)^{n-1}(t-\ln(2t))^2 = 3(t-\ln(2t))^2. \nonumber \\
\end{eqnarray}
The last factor is bounded from above and below for $t \in [1/2,1]$, which proves the lemma. 
\noindent
\end{proof}
Similar behaviour is observed for ordinary critical percolation (see~\cite{Grim}, Section 10.1).  
\label{crittree} 
\end{section}

\begin{section}{Dependencies}
Any two sites $v$ and $w$ in the tree, are connected by a unique path. 
The dependence between the colours of sites $v$ and $w$ at time $t$, decays exponentially in $\pi(v,w)$, the number of sites between $v$ and $w$. This is also true for the modified frozen percolation process  on the line $\mathbb{Z}$. 
\begin{lemma}
There exists a constant $A \in (0,1)$ such that for all sites $v$ and $w$ and all $t \in [0,1]$, 
and for any $\mathcal{C}_k(t) \in \{\mathcal{W}(t), \mathcal{G}(t), \mathcal{R}(t)\}, k = 1,2$,
\begin{equation}
\Bigl| \mathcal{P}(v \in \mathcal{C}_1(t), w \in \mathcal{C}_2(t)) - \mathcal{P}(v \in \mathcal{C}_1(t))\mathcal{P}(w \in \mathcal{C}_2(t)) \Bigr| \leq 5A^{\lfloor\frac{\pi(v,w)}{12}\rfloor}.
\label{indeptree}
\end{equation}
\label{treedep}
\end{lemma}
\begin{proof}
For $t \leq 1/2$, sites can only be white or green and the transition from white to green happens independent of the other sites. 
In fact, the left hand side of \eqref{indeptree} equals zero. So suppose $t > 1/2$ and let $v = i_1, i_2, \ldots, i_n = w$ denote the unique path 
between $v$ and $w$. We may assume that $n \geq 12$, otherwise there is nothing to prove. Define the following events:
\begin{eqnarray}
E_1 &:=& \{ \exists j \in (3n/4,n) : \min \{U_{i_j}, U_{i_{j+2}}\} > Z_{i_{j+1}} \} \nonumber \\
E_2 &:=& \{ \exists k \in (0,n/4) : \min \{ U_{i_k}, U_{i_{k-2}} \} > Z_{i_{k-1}} \} \nonumber 
\end{eqnarray}
If $\min \{U_{i_j}, U_{i_{j+2}}\} > Z_{i_{j+1}}$ for some $j$, it follows that $Z_{i_{j+1}} = Y_{i_{j+1} \rightarrow a_{j+1}}$, where $a_{j+1}$ is the child of $i_{j+1}$ not lying on the path from $v$ to $w$. 
So if $E_1$ occurs for some $j$, it only depends on $U_{i_j}, U_{i_{j+1}}, U_{i_{j+2}}$ and the $U$-variables in the subtree leaving from $a_{j+1}$. Fix $j \in (3n/4,n)$. Taking disjoint sets of three consecutive sites on the path from $v$ to $w$, we arrive at  
\begin{equation}
\mathcal{P}(E_1^c) \leq \bigl(1 - \mathcal{P}(Y_{i_{j+1}\rightarrow a_{j+1}} < \min \{U_{i_j}, U_{i_{j+2}} \}\bigr)^{\lfloor\frac{n}{12}\rfloor},
\label{boundec1} 
\end{equation}
by automorphism invariance. 
A similar reasoning holds for $E_2$. 
Further, on $E_1 \cap E_2$, the colours of $v$ and $w$ are independent: suppose we cut the tree in two parts by removing any edge $(i_l,i_{l+1})$ 
with $n/4 < l < 3n/4-1$ on the path from $v$ to $w$. The occurrence of $E_1$ and $E_2$ ensure that the $U$-variables on the part containing $v$ ($w$) determine the colour of $v$ ($w$, respectively) uniquely. If $E_1$ occurs for some $j$, at least one of the sites $i_j$ and $i_{j+2}$ is at rest at time $t$, or site $i_{j+1}$ is frozen without $i_j$ or $i_{j+2}$ at time $t$. In both cases, the configuration on $i_j,i_{j+1},i_{j+2}$ prevents anything happening on the part of the tree containing $v$ to influence the colour of $w$. We compute 
\begin{eqnarray}
&& \hspace{-30pt}
\mathcal{P}(v \in \mathcal{C}_1(t), w \in \mathcal{C}_2(t)) \nonumber \\
&=& \mathcal{P}(v \in \mathcal{C}_1(t), w \in \mathcal{C}_2(t), E_1, E_2) + \mathcal{P}(v \in \mathcal{C}_1(t), w \in \mathcal{C}_2(t), E_1^c \cup  E_2^c ) \nonumber \\
&=& \mathcal{P}(v \in \mathcal{C}_1(t) \cap E_2) \mathcal{P}(w \in \mathcal{C}_2(t) \cap E_1) +  \mathcal{P}(v \in \mathcal{C}_1(t), w \in \mathcal{C}_2(t), E_1^c \cup  E_2^c ) \nonumber \\
&=& \bigl[ \mathcal{P}(v \in \mathcal{C}_1(t)) - \mathcal{P}(v \in \mathcal{C}_1(t) \cap E_2^c) \bigr]\bigl[\mathcal{P}(w \in \mathcal{C}_2(t)) -  \mathcal{P}( w \in \mathcal{C}_2(t) \cap E_1^c)\bigr]  \nonumber \\
&&+ \ \   \mathcal{P}(v \in \mathcal{C}_1(t), w \in \mathcal{C}_2(t), E_1^c \cup  E_2^c ). \label{tocopy}
\end{eqnarray}
From \eqref{tocopy} we obtain
\begin{equation}
\Bigl| \mathcal{P}(v \in \mathcal{C}_1(t), w \in \mathcal{C}_2(t)) - \mathcal{P}(v \in \mathcal{C}_1(t))\mathcal{P}(w \in \mathcal{C}_2(t)) \Bigr| \leq 5\mathcal{P}(E_1^c). 
\end{equation}

The only thing we need to show now, is that the probability on the left side of \eqref{boundec1} is bounded away from $1$. 
We use that by \eqref{defF}, 
$$\mathcal{P}(Y_{i_{j+1}\rightarrow a_{j+1}} < \min\{U_{i_j},U_{i_{j+2}}\}) = 2\int_{1/2}^1 (1-s)\ln(2s) ds > 0.$$  This allows us to define the required constant $A$.
\end{proof}   
  
\begin{rem}
One may wonder why we took the site version of the model, whereas Aldous' model is stated in terms of edges. The first reason is that we studied the frozen percolation model on $\mathbb{Z}$ using the site version and we wanted to keep the models consistent. The second reason is that we wanted to see whether a perturbation in the dynamics would cause the critical behaviour to disappear. The third (and in practice maybe the most important) reason is that the edge version introduces dependencies, which make it harder (if not impossible) to construct the model.

 Up to a certain point, we could just as well have taken the edge version of the process. Then, we replace each edge $e$ of the tree by two directed edges $\re$ and $\el$. We attach an independent $U[0,1]$ variables $U_e = U_{\re} = U_{\el}$ to each edge. Intuitively, we need variables Y$_{\re}$ and $Y_{\el}$ that represents the time a directed edge freezes in the directed subtree, where the direction is inherited from $\re$ or $\el$. We can define a relation similar to \eqref{defy} which $Y_{\re}$ and $Y_{\el}$ have to satisfy (recall the definition of $\Phi$ in \eqref{defPhi}): 
\begin{equation}
Y_{\re} = \Phi(\min(Y_{\stackrel{\rightarrow}{e_1}},Y_{\stackrel{\rightarrow}{e_2}}),\max(Y_{\stackrel{\rightarrow}{e_1}},Y_{\stackrel{\rightarrow}{e_2}}), U_e),
\label{defe}
\end{equation}
where $\stackrel{\rightarrow}{e_1}$ and $\stackrel{\rightarrow}{e_2}$ are the children of the directed edge $\re$. We follow the heuristic arguments from Section \ref{heuristics}, using \eqref{defe}. 
 The distribution $F(t) = \ln(2t)$ for $t\geq 1/2$ that drops out (compare Lemma \ref{consistency}) for the edge version is the same as for the site version. The analogue of Lemma \ref{existencelemma} can be proved similarly. However, problems arise because the edge version of frozen percolation carries more dependencies that the site version. 
Particularly, the analogue of $Z_i$ (see equation \eqref{defZ}) for the edge version of the model is not easily defined. 
\label{edge}
\end{rem}
\begin{opques}
Is it possible to construct the edge version of the frozen percolation process?
\end{opques}


\end{section}
\label{froztree}
\bibliographystyle{alpha}
\bibliography{thesis}

\end{document}